\documentclass[reqno]{amsart}

\usepackage{latexsym,amssymb,amsthm,amsmath,amscd}

\theoremstyle{plain}
\newtheorem{theorem}{Theorem}[section]

\newtheorem{lemma}{Lemma}[section]
\newtheorem{proposition}{Proposition}[section]

\newtheorem{corollary}{Corollary}[section]

\newtheorem{definition}{Definition}[section]
\newtheorem{example}{Example}[section]

\theoremstyle{remark}
\newtheorem{remark}{Remark}[section]

\numberwithin{equation}{section}

\def\<{\left < }
\def\>{\right >}
\def\({\left ( }
\def\){\right )}
\def\e{\eqref}

\def\e{\eqref}

\begin{document}
\title[Geometry of $CRS$ bi-warped product submanifolds]{Geometry of $CRS$ bi-warped product submanifolds in Sasakian and cosymplectic manifolds}
\author[B.-Y. Chen]{Bang-Yen Chen}
\address{Department of Mathematics, Michigan State University, East Lansing, Michigan,
48824--1027, U.S.A.}
\email{chenb@msu.edu}
\author[S. Uddin]{Siraj Uddin}
\address{Department of Mathematics, Faculty of Science, King Abdulaziz University,Jeddah 21589, Saudi Arabia}
\email{siraj.ch@gmail.com}

\author[A. Alghanemi]{Azeb Alghanemi}
\address{Department of Mathematics, Faculty of Science, King Abdulaziz University, Jeddah 21589, Saudi Arabia}
\email{aalghanemi@kau.edu.sa}

\author[A. Al-Jedani]{Awatif Al-Jedani}
\address{Department of Mathematics, Faculty of Science, University of Jeddah, 23890 Jeddah, Saudi Arabia}
\email{amaljedani@uj.edu.sa}
\author[I. Mihai]{Ion Mihai}
\address{Faculty of Mathematics and Computer Science, University of Bucharest, Str. Academiei 14, 010014 Bucharest, Romania}
\email{imihai@fmi.unibuc.ro}

\begin{abstract} In this paper, we prove that there are no proper $CRS$ bi-warped product submanifolds other than contact CR-biwarped products in Sasakian manifolds. On the other hand, we prove that if $M$ is a $CRS$ bi-warped product of the form $M=N_T \times_{f_1}N^{n_{1}}_\perp\times_{f_2} N^{n_{2}}_\theta$ in a cosymplectic manifold $\widetilde M$, then its second fundamental form $h$ satisfies the inequality: $$\|h\|^2\geq 2n_1\|\nabla(\ln f_1)\|^2+2n_2(1+2\cot^2\theta)\|\nabla(\ln f_2)\|^2,$$ where $N_T,\, N^{n_{1}}_\perp$ and $N^{n_{2}}_\theta$ are invariant, anti-invariant and proper pointwise slant submanifolds of $\widetilde M$, respectively,  and $\nabla(\ln f_1)$ and $\nabla(\ln f_2)$ denote the gradients of $\ln f_{1}$ and $\ln f_{2}$, respectively. Several applications of this inequality are given. At the end, we provide a non-trivial example of bi-warped products satisfying the equality case.
\end{abstract}
\maketitle

\section{Introduction}

Let $B,N_{1}$ and $N_{2}$ be Riemannian manifolds equipped with metrics $g_B,g_1$ and $g_{2}$, respectively, and let $f_{1},f_{2}$ be two positive functions on $B$. Then the {\it bi-warped product} $B\times_{f_{1}} N_{1}\times_{f_{2}} N_{2}$ is the product manifold $B\times N_{1}\times N_{2}$ equipped with the bi-warped product metric
\begin{equation}\label{E:warped} g=g_B+f_{1}^2 g_1+f_{2}^2 g_2.\end{equation}
Bi-warped products extend the ``ordinary'' warped products and they form  
an important class of multiply warped products (see, e.g., \cite{Awatif19,C6,CD,SC1,SC,Unal}).

A submanifold $M$ of a Kaehler manifold $(\widetilde M, \widetilde g,J)$ is called a {\it $CR$ submanifold}  \cite{bejancu} if there exist a holomorphic distribution $\mathfrak D^{T}$ and a totally real distribution $\mathfrak D^\perp$ on $M$ such that $TM=\mathfrak D^{T}\oplus \mathfrak D^\perp$, where $TM$ denotes the tangent bundle of $M$. 

The study of warped products from submanifold theory points of view was initiated by
the first author in \cite{C3,C4,C02,C03,C04}. In particular, he proved in \cite{C3} that there do not exist warped products of the form $N_{\perp}\times_{f} N_{T}$ in any Kaehler manifold, where $N_{T}$ is a holomorphic submanifold and $N_{\perp}$ is a totally real submanifold. On the other hand, he proved that, for any  warped product of the form $M=N_{T}\times_f N^{n_{1}}_{\perp}$ in any Kaehler manifold,  the second fundamental form $h$ of $M$ satisfies
$$||h||^2\geq 2n_1 ||\nabla(\ln f)||^2,$$ where $n_{1}=\dim\,N^{n_{1}}_{\perp}$ and $\nabla (\ln f)$ is the gradient of $\ln f$. He simply  called such warped product submanifolds $N_{T}\times_f N^{n_{1}}_{\perp}$ to be {\it $CR$-warped products}.

In this paper, we investigate an important class of bi-warped product submanifolds in a contact metric manifold $\widetilde M$; namely, $CRS$ bi-warped products of the form $N_T\times_{f_1}N^{n_{1}}_\perp\times_{f_2}N^{n_{2}}_\theta$, where $N_T, N^{n_{1}}_\perp$ and $N^{n_{2}}_\theta$ are invariant, anti-invariant and pointwise slant submanifolds of $\widetilde M$ and $f_{1},f_{2}$ are warping functions.
 We simply call such a submanifold a {\it $CRS$ bi-warped product}. 

In this paper, first we prove that there are no proper $CRS$ bi-warped products other than contact $CR$-biwarped products in any Sasakian manifold. Then we investigate $CRS$ bi-warped products $M=N_T\times_{f_1}N^{n_{1}}_\perp\times_{f_2}N^{n_{2}}_\theta$ in cosymplectic manifolds. The main result of this paper states that, for such a $CRS$ bi-warped product in a cosymplectic manifold, its second fundamental form $h$ satisfies
 $$\|h\|^2\geq 2n_1\|\nabla(\ln f_1)\|^2+2n_2(1+2\cot^2\theta)\|\nabla(\ln f_2)\|^2,$$
where $n_1=\dim\,N^{n_{1}}_\perp,\,n_2=\dim\, N^{n_{2}}_\theta$ and $\theta$ is the slant function of $N^{n_{2}}_\theta$.
In this paper, we also study the equality case of the inequality. Several applications of our main result are also given. Furthermore, we provide a non-trivial example of submanifold which verify the equality case of the inequality.


\section{Preliminaries}\label{S2}

A $(2m+1)$-dimensional Riemannian  manifold $(\widetilde M^{2m+1}, g)$ is said to be a {\it{contact metric manifold}} if it admits a $(1, 1)$ tensor field $\varphi$ on its tangent bundle $T\widetilde M^{2m+1}$, a vector field ${\mathfrak{\xi}}$ (i.e., the structure  vector field) and a $1$-form $\eta$ satisfying:
\begin{equation}\begin{aligned}\label{struc}
&\varphi^2=-Id+\eta\otimes{\mathfrak{\xi}},\,\,\,\eta({\mathfrak{\xi}})=1,\,\,\,\varphi(\xi)=0,\,\,\,\eta\circ\varphi=0,\\
&g(\varphi X, \varphi Y)=g(X, Y)-\eta(X)\eta(Y),\,\,\,\eta(X)=g(X, {\mathfrak{\xi}}).
\end{aligned}\end{equation}
for any vector fields $X, Y$ on $\widetilde M^{2m+1}$, where $\widetilde\nabla$ denotes the Riemannian connection with respect to the Riemannian metric $g$ (see \cite{Bl}). An almost contact metric manifold $(\widetilde M, g, \varphi, \xi, \eta)$ is said to be {\textit{normal}} if the tensor field $N_{\varphi}=[\varphi, \varphi]+2d\eta\otimes\xi$ vanishes identically, where $[\varphi, \varphi]$ denotes the Nijenhuis tensor of $\varphi$. 

A normal contact metric manifold is said to be a {\textit{Sasakian manifold}}. An almost contact metric manifold is Sasakian if and only if
\begin{align}\label{sasaki}
(\widetilde\nabla_X\varphi)Y=g(X, Y){\mathfrak{\xi}}-\eta(Y)X,\,\,\,\widetilde\nabla_X{\mathfrak{\xi}}=-\varphi X,
\end{align}
for any $X, Y\in\Gamma(T\widetilde M)$, where $\Gamma(T\widetilde M)$ is the Lie algebra of vector fields on $\widetilde M$.
An almost contact metric manifold is called {\em{almost cosymplectic}}  if $d\eta=0$ and $d\varphi=0$ according to D. E. Blair in \cite{Bl}. In particular, an almost cosymplectic manifold is called {\em{cosymplectic}} if it satisfies
\begin{align}\label{cosy}
\widetilde\nabla\varphi=0,\quad
\widetilde\nabla\xi=0.
\end{align}

Let $M^n$ be an $n$-dimensional submanifold of an almost contact metric manifold  $\widetilde M^{2m+1}$ such that the structure vector field ${\mathfrak{\xi}}$ tangent to $M^n$ with induced metric $g$. Let $\Gamma(T^\perp M^{n})$ denote the space of all vector fields normal to $M^n$. Then the Gauss and Weingarten formulas are given respectively by (see, e,g., \cite{book11,book17, Ya})
\begin{equation}\begin{aligned}\label{gauss}
\widetilde\nabla_XY=\nabla_XY+h(X, Y), \quad
\widetilde\nabla_XN=-A_NX+\nabla^{\perp}_XN, 
\end{aligned}\end{equation}
 for  $X, Y\in\Gamma(TM^n)$ and $N\in\Gamma(T^{\perp}M^{n})$, where $\nabla$ and $\nabla^{\perp}$ denote the connections on the tangent and normal bundles of $M^n$, respectively, and $A$ the shape operator of $M^n$. It is well known that the second fundamental form $h$ and the shape operator $A$ are related by $g(h(X, Y), N)=g(A_NX, Y).$

Let $x\in M^n$ and $e_1,\cdots, e_n, e_{n+1},\cdots, e_{2m+1}$ are orthonormal basis of the tangent space $\widetilde M^{2m+1}$ such that restricted to $M^n$,  $e_1,\cdots, e_n$ are tangent to $M^n$ and hence $e_{n+1},\cdots, e_{2m+1}$ are normal to $M^n$. Let $\{h^r_{ij}\}\,(1\leq i, j\leq n;\,n+1\leq r\leq 2m+1)$ denote the coefficients of $h$ with respect to the local frame field. Then
\begin{align} \label{funda}
h_{ij}^r=g(h(e_i, e_j), e_r)=g(A_{e_r}e_i, e_j),\,\,\,\|h\|^2=\sum_{i,j=1}^ng(h(e_i, e_j), h(e_i, e_j)).\end{align}
The mean curvature vector $H$ is defined by $H=\frac{1}{n} \sum_{i=1}^{n} h(e_i, e_i)$. A submanifold $M$ is called {\it totally geodesic} if $h=0$; {\it totally umbilical} if $h(X, Y)=g(X, Y)H$; and {\it minimal} if $H=0$.

For submanifolds $M$ of an almost contact metric manifold $(\widetilde M^{2m+1},g,\varphi,\xi,\eta)$, there are three important classes of submanifolds depending on the action of $\varphi$ on $TM$. 
\vskip.07in

\noindent (a) A submanifold $M^{2n+1}$ tangent to ${\mathfrak{\xi}}$ is called {\textit{invariant}}  if $\varphi$ preserves each tangent space of $M^{2n+1}$, i.e., $\varphi(T_xM^{2n+1})\subseteq T_xM^{2n+1}, \forall x\in M^{2n+1}$.
\vskip.05in

\noindent (b)  A submanifold $M^{n}$ is called {\textit{anti-invariant}}  if $\varphi$ maps any tangent space of $M^{n}$ into the corresponding normal space, i.e., $\varphi(T_xM^{n})\subset T^\perp_xM^{n}, \forall x\in M^{n}$.
\vskip.05in

\noindent (c)  A submanifold $M^{n}$ is called {\textit{pointwise slant}}  \cite{Etayo,CG} if for any nonzero vector $X\in T_xM^n\,(x\in M)$,  the angle $\theta(X)$ between $\varphi X$ and $T_xM^n$ is independent of the choice of $X\in T_xM$. In this case, $\theta$ defines a function on $M^n$, called the {\textit{slant function}}. In particular, if the slant function $\theta$ is a global constant on $M^{n}$, then $M^{n}$ is said to be a {\it slant submanifold} or a {\it $\theta$-slant submanifold}. 
\vskip.07in

Anti-invariant submanifolds are slant submanifolds with slant function $\theta=\frac{\pi}{2}$. A pointwise slant submanifold is called {\it proper} if $0<\theta <\frac{\pi}{2}$. See \cite{Park,UK} for  non-trivial examples of pointwise slant submanifolds.

For any $X\in\Gamma(TM^n)$ and $N\in \Gamma(T^{\perp}M^n)$, we put $$\varphi X = T X + FX,\quad \varphi N = tN+fN,$$ where $TX$ (resp., $tN$) is the tangential  component, and $FX$ (resp., $fN$) is the  normal component of $\varphi X$ (resp. of $\varphi N$).

We recall the following useful characterization from \cite{UK}.

\begin{proposition}\label{P:2.1} Let $M^n$ be a submanifold of an almost contact metric manifold $\widetilde M^{2m+1}$ with $\mathfrak{\xi}\in\Gamma(TM^n)$. Then $M^n$ is pointwise slant if and only if 
\begin{align}\label{charac} T^2=\cos^2\theta\left(-I+\eta\otimes\mathfrak{\xi}\right),\end{align}
where $\theta$ is the slant function and  $I$ denotes the identity map on $TM^{n}$.
\end{proposition}

Following relations are straightforward consequence of \eqref{charac}
\begin{align}
&g(TX,TY)=\cos^2\theta[g(X,Y)-\eta(X)\eta(Y)], \label{slant1}\\
&g(FX, FY) =\sin^2\theta[g(X, Y)-\eta(X)\eta(Y)],\label{slant2}
\end{align}
for  vector fields $X,Y\in \Gamma(M^n)$.
Also, for pointwise slant submanifolds, we have 
\begin{align}
 \label{slant3}
tFX=\sin^2\theta\left(-X+\eta(X)\mathfrak{\xi}\right),~~~~fFX=-FTX,\quad X\in\Gamma(TM^n).
\end{align}

Next, we give the following useful result for later use.

\begin{proposition}\label{P:2.1CS} Every pointwise slant submanifold $M$ of a Sasakian manifold $\widetilde M^{2m+1}$ is anti-invariant if and only if the structure vector field $\xi$ is normal to $M$.
\end{proposition}
\begin{proof}
From \eqref{sasaki} and \eqref{gauss}, we have
\begin{align}
 \label{cs1}
 A_\mathfrak{\xi}X=TX,\quad
 \nabla^\perp_X\xi=-FX,
 \end{align}
 for any $X\in\Gamma(TM)$. Taking the inner product in \eqref{cs1} with $Y\in\Gamma(TM)$, we get
 \begin{align}
 \label{cs2}
 g(h(X, Y), \mathfrak{\xi})=g(TX, Y).
 \end{align}
 From \eqref{cs2}, we find $\cos^2\theta\;g(X, Y)=0$ via polarization identity. Since, the metric $g$ is Riemannian,  we obtain $\theta=\frac{\pi}{2}$. The converse is trivial. 
\end{proof}

\section{Some basic results on $CRS$ bi-warped products}\label{S3}

Let $M=B\times_{f_1}N$ with $N=N_1\times_{f_2}N_2$ be a bi-warped product submanifold of a Riemannian manifold $\widetilde M$. Put ${\mathfrak{D}}=TB$. Then we have
\begin{align}\label{warped} \nabla_ZX= \nabla_XZ=\sum_{i=1}^2X(\ln f_i) Z^i \end{align}
for any $X\in \Gamma({\mathfrak{D}})$ and $Z\in\Gamma(TN)$, where  $Z^i$ denotes the $N_i$-component of $Z$, and $\nabla$ is the Levi-Civita connection on $M$ (see, e.g., \cite{UFH}).

\begin{remark}\label{R:3.1} It is well known that, for a warped product manifold $M=B\times_fF$ , $B$ is totally geodesic and $F$ is totally umbilical in $M$ (see, e.g., \cite{C3,oneill}).
\end{remark}

Now, we give the following definition of $CRS$ bi-warped product submanifolds.

\begin{definition}\label{D:3.1} {\rm A {\it $CRS$ bi-warped product} of an almost contact metric manifold $(\widetilde M^{2m+1},g,\varphi,\xi,\eta)$ is a bi-warped product $M^n=N_{T}^{2p+1}\times_{f_1}N^{n_1}_\perp\times_{f_2}N^{n_2}_\theta$ in $\widetilde M^{2m+1}$, where $N_{T}^{2p+1}$ is an invariant submanifold tangent to $\xi$, $N^{n_1}_\perp$ is an anti-invariant submanifold,  $N^{n_2}_\theta$ is a pointwise slant submanifold of $\widetilde M^{2m+1}$, and $f_{1}$ and $f_{2}$ are two positive smooth functions on $N_{T}^{2p+1}$.} \end{definition}

\begin{remark}\label{R:3.2} It follows from Remark \ref{R:3.1} that $N_{T}, N_{T}\times_{f_1}N^{n_1}_\perp$ and $N_{T}\times_{f_2}N^{n_2}_\theta$ are totally geodesic in the $CRS$ bi-warped product $M^n=N_{T}\times_{f_1}N^{n_1}_\perp\times_{f_2}N^{n_2}_\theta$.
\end{remark}

\begin{remark}\label{R:3.3} A $CRS$ bi-warped product is a {\em{contact CR-biwarped product}} if it is of the form $M^n=N_{T}^{2p+1}\times_{f_1}N^{n_1}_\perp\times_{f_2}N^{n_2}_\perp$, where $N_{T}^{2p+1}$ is an invariant submanifold tangent to $\xi$, $N^{n_1}_\perp$ and $N^{n_2}_\perp$ are anti-invariant submanifolds of $\widetilde M^{2m+1}$.
\end{remark}

Let ${\mathfrak{D}^{T}},\,\,{\mathfrak{D}}^\perp$ and ${\mathfrak{D}}^\theta$ denote the tangent bundles of $N_T,\,\,N_\perp$ and $N_\theta$, respectively. Then the tangent and normal spaces of $M^n$ are given respectively by
 \begin{align}\label{tan-nor}
 & TM^n={\mathfrak{D}^{T}}\oplus{\mathfrak{D}}^\perp\oplus{\mathfrak{D}}^\theta, \quad
T^\perp M^n=\varphi{\mathfrak{D}}^\perp\oplus F{\mathfrak{D}}^\theta\oplus\mu,
  \end{align}
where  $\mu$ is the invariant subbundle of $T^\perp M^n$. Clearly, ${\rm Span}\{\xi\}\subset {\mathfrak{D}^{T}}$ by definition.

The following result provides the non-existence of pointwise semi-slant warped products in Sasakian manifolds.

\begin{proposition}\label{P:3.1} There do not exist $CRS$ bi-warped products in a Sasakian manifold other than contact $CR$-biwarped products.
\end{proposition}
\begin{proof} Let $M=N_{T}\times_{f_1}N_\perp\times_{f_2}N_\theta$ be a $CRS$ bi-warped product in a Sasakian manifold $\widetilde M$. 
 If $V\in \Gamma({\mathfrak{D}}^\theta)$, it follows from \e{struc}, \e{sasaki}, \e{gauss} and \e{warped} that
 \begin{align}\label{3.04}\xi(\ln f_{2})V+ h(\xi,V)=\nabla_{V}\xi +h(\xi,V)=\widetilde \nabla_{V}\xi=-\varphi V=- TV -FV.\end{align}
By comparing the tangential components of \e{3.04}, we find $\xi(\ln f_{2})V=-TV$, which implies $\xi(\ln f_{2})=0$ and $TV=0$. Hence, we find $\theta=\frac{\pi}{2}$. Consequently, $M$ is a contact $CR$-biwarped product.
\end{proof}

The following result is an immediate consequence of Proposition \ref{P:3.1}.
\begin{corollary}\label{C:3.1CS} There do not exist any proper pointwise semi-slant warped product submanifold $M=N_T\times_fN_\theta$ in a Sasakian manifold.
\end{corollary}

\begin{remark} One may also prove that there are no proper bi-warped products of the form $M=N_{T}\times_{f_{1}} N_{\perp}\times_{f_{2}} N_{\theta}$ in a Sasakian manifold such that the characteristic vector field $\xi$ is tangent to either $N_{\perp}$ or to $N_{\theta}$, or $\xi$ is normal to $M$.
\end{remark}

In the next section, we have found that if we replace the ambient  spaces from  Sasakian to cosymplectic, then such warped products  $N_T\times_{f_1}N_\perp\times_{f_2}N_\theta$ do exist.

\section{$CRS$ bi-warped products of cosymplectic manifolds}\label{S4}

Proposition \ref{P:3.1} shows that there are no proper $CRS$ bi-warped products in any Sasakian manifold. In this section, we study  $CSR$ bi-warped products in cosymplectic manifolds which is quite different from Sasakian case. 

From now on, assume that  $M=N_{T}\times_{f_1}N_\perp\times_{f_2}N_\theta$ is a $CRS$ bi-warped product of a cosymplectic manifold $\widetilde M$ whose structure vector field $\xi$ is tangent to $N_T$. 

Now, we give the following useful lemmas for later use.

\begin{lemma}\label{L:4.1} Let $M^n=N_T\times_{f_1}N_\perp\times_{f_2}N_\theta$ be a $CRS$ bi-warped product in a cosymplectic manifold $\widetilde M^{2m+1}$. Then $\xi(\ln f_1)=\xi(\ln f_2)=0$ and $h(X, \xi)=0$ for every $X\in \Gamma(N_{\perp}\times N_{\theta})$.
\end{lemma}
\begin{proof} It follows from  \e{gauss} and \e{cosy} that 
$\nabla_{X}\xi+h(X,\xi)=\widetilde \nabla_X {\xi}=0$ for every $X\in \Gamma(N_{\perp}\times N_{\theta})$, which implies $\nabla_{X}\xi=0$ and $h(X,\xi)=0$. Now, after combining $\nabla_{X}\xi=0$  with \e{warped}, we obtain $\xi(\ln f_1)=\xi(\ln f_2)=0$.
\end{proof}

\begin{lemma}\label{L:4.2} Let $M^n=N_T\times_{f_1}N_\perp\times_{f_2}N_\theta$ be a $CRS$ bi-warped product in a cosymplectic manifold $\widetilde M^{2m+1}$. Then we have
\begin{enumerate}
\item[(i)] $g(h(X, Y), \varphi Z)=0=g(h(X, Y),FV)$;
\item [(ii)] $g(h(X, Z), \varphi W)=-\varphi X(\ln f_1)\,g(Z, W)$;
\item[(iii)] $g(h(X,U), FV)=-\varphi X(\ln f_2)\,g(U, V)-X(\ln f_2)\,g(U, TV)$.
\end{enumerate}
for $X,Y\in \Gamma({\mathfrak{D}^{T}})$, $Z,W\in \Gamma({\mathfrak{D}}^\perp)$ and  $U,V\in \Gamma({\mathfrak{D}}^\theta)$.
\end{lemma}
\begin{proof}
For any $X, Y\in \Gamma({\mathfrak{D}^{T}})$ and $Z\in \Gamma({\mathfrak{D}}^\perp)$, we find
\begin{align*}
g(h(X, Y), \varphi Z)\, &=g(\widetilde\nabla_XY, \varphi Z)=-g(\widetilde\nabla_X\varphi Y, Z)=-g(\nabla_X\varphi Y, Z).
\end{align*}
Using the fact that $N_T$ is totally geodesic in $M^n$ (see Remark \ref{R:3.1}) with the orthogonality of vector fields, we find $g(h(X, Y), \varphi Z)=0$, which is the first equality of (i). Similarly, we have the second equality of (i). 
For (ii), we have
\begin{align*}
g(h(X, Z), \varphi W)\, &=g(\widetilde\nabla_ZX, \varphi W)=g((\widetilde\nabla_Z\varphi)X, W)-g((\widetilde\nabla_Z\varphi X, W).
\end{align*}
Using \eqref{cosy}, Remark \ref{R:3.2}, and \eqref{warped} with the orthogonality of vector fields, we get (ii). Similarly, we obtain (iii).  \end{proof}

Replacing $X$ by $\varphi X$  in Lemma \ref{L:4.2}(ii) and using \eqref{struc} and  Lemma \ref{L:4.1}, we find
 \begin{align}\label{4.1}
g(h(\varphi X, Z), \varphi W)=X(\ln f_1)g(Z, W).
\end{align}

Similarly, if we replace $X$ by $\varphi X$, $U$ by $TU$, and $V$ by $TV$ in Lemma \ref{L:4.2}(iii), we obtain the following useful relations:
 \begin{align}\label{4.2}
&g(h(\varphi X, U), FV)=X(\ln f_2)\,g(U, V)-\varphi X(\ln f_2)\,g(U, TV),
\\ \label{4.3}
&g(h(X,TU), FV)=\varphi X(\ln f_2)\,g(U, TV)-\cos^2\theta\,X(\ln f_2)\,g(U, V),
\\ &\label{4.4} g(h(X,U), FTV)=\cos^2\theta\,X(\ln f_2)\,g(U, V)-\varphi X(\ln f_2)\,g(U, TV),
\\&\label{4.5} g(h(\varphi X,TU), FV)=-X(\ln f_2)\,g(U, TV)-\cos^2\theta\,\varphi X(\ln f_2)\,g(U, V),
\\&\label{4.6} g(h(\varphi X, U), FTV)=X(\ln f_2)\,g(U, TV)+\cos^2\theta\,\varphi X(\ln f_2)\,g(U, V),
\\ \label{4.7} &g(h(X, TU), FTV)=-\cos^2\theta\,X(\ln f_2)\,g(U, TV)-\cos^2\theta\,\varphi X(\ln f_2)\,g(U, V),
\\&\label{4.8} g(h(\varphi X, TU), FTV)=\cos^2\theta\,X(\ln f_2)\,g(U, V)-\cos^2\theta\,\varphi X(\ln f_2)\,g(U, TV).
\end{align}

\begin{lemma}\label{L:4.3} Let $M^n=N_T\times_{f_1}N_\perp\times_{f_2}N_\theta$ be a $CRS$ bi-warped product in a cosymplectic manifold $\widetilde M^{2m+1}$. Then \begin{align}\label{4.9}
g(h(X, Z),FV)=g(h(X, V), \varphi Z)=0,\end{align}
for any $X\in \Gamma({\mathfrak{D}^{T}})$, $Z\in \Gamma({\mathfrak{D}}^\perp)$ and $V\in \Gamma({\mathfrak{D}}^\theta)$.
\end{lemma}
 \begin{proof} For any $X\in \Gamma({\mathfrak{D}^{T}})$,\,\, $Z\in \Gamma({\mathfrak{D}}^\perp)$ and $V\in \Gamma({\mathfrak{D}}^\theta)$, we have
\begin{equation*} \begin{aligned}
g(h(X, Z), FV)\;&=g(\widetilde{\nabla}_ZX-\nabla_ZX, \varphi V-TV)\\&=g((\widetilde{\nabla}_{Z}\varphi)X, V) -g(\widetilde\nabla_{Z}\varphi X, V)-X(\ln f_1)\,g(Z, TV).
\end{aligned}\end{equation*}
Then we find the first equality of \eqref{4.9} from \eqref{cosy}, \e{warped} and the orthogonality of vector fields. 
For the second equality, first we have
 \begin{align*}
g(h(X, V), \varphi Z)&=g(\widetilde{\nabla}_VX, \varphi Z)=-g(\widetilde\nabla_V\varphi X, Z)=-\varphi X(\ln f_2)\,g(Z, V)=0,
\end{align*}
which gives the second equality of \eqref{4.9}. 
 \end{proof}

\begin{proposition}\label{P:4.1}
Let $M^n=N_T\times_{f_1}N_\perp\times_{f_2}N_\theta$ be a $CRS$ bi-warped product in a cosymplectic manifold $\widetilde M^{2m+1}$ such that $h({\mathfrak{D}^{T}}, {\mathfrak{D}}^\perp)\perp\mu$. Then the warping function $f_1$ is constant if and only if $M^n$ is ${\mathfrak{D}^{T}}\oplus{\mathfrak{D}}^\perp$-mixed totally geodesic in $\widetilde M^{2m+1}$, i.e., $h({\mathfrak{D}^{T}}, {\mathfrak{D}}^\perp)=\{0\}$.
\end{proposition}
 \begin{proof} If $M^n$ is ${\mathfrak{D}^{T}}\oplus{\mathfrak{D}}^\perp$-mixed totally geodesic,  \eqref{4.2} implies that $f_1$ is constant.

Conversely, if $f_1$ is constant, then it follows from Lemma \ref{L:4.3} that
 \begin{align}\label{4.10}
g(h(\mathfrak{D}^{T}, \mathfrak{D}^\perp),F\mathfrak{D}^\theta)=0.
\end{align}
On the other hand, we find from \eqref{4.2} that
$g(h(\mathfrak{D}^{T}, \mathfrak{D}^\perp), \varphi \mathfrak{D}^\perp)=0.$
Also, from the hypothesis of the proposition we have
$g(h(\mathfrak{D}^{T}, \mathfrak{D}^\perp), \mu)=0.$
Thus, we get $g(h(\mathfrak{D}^{T}, \mathfrak{D}^\perp), N)=0$ for any $N\in\Gamma(T^\perp M^n)$. Therefore, $M^n$ is $\mathfrak{D}^{T}\oplus\mathfrak{D}^\perp$-mixed totally geodesic. 
 \end{proof}

\begin{proposition}\label{P:4.2}
Let $M^n=N_T\times_{f_1}N_\perp\times_{f_2}N_\theta$ be a $CRS$ bi-warped product in a cosymplectic manifold $\widetilde M^{2m+1}$ such that $h({\mathfrak{D}^{T}}, {\mathfrak{D}}^\theta)\perp\mu$. Then $f_2$ is constant on $M$ if and only if $M^n$ is ${\mathfrak{D}^{T}}\oplus{\mathfrak{D}}^\theta$-mixed totally geodesic in $\widetilde M^{2m+1}$.
\end{proposition}
 \begin{proof} If $M^n$ is ${\mathfrak{D}^{T}}\oplus{\mathfrak{D}}^\theta$-mixed totally geodesic, then \eqref{4.2} yields
 \begin{align}\label{4.13}
X(\ln f_2)\,g(U, V)=\varphi X(\ln f_2)\,g(U, TV).
\end{align}
Also, from \eqref{4.3} we find
$\varphi X(\ln f_2)\,g(U, TV)=\cos^2\theta\,X(\ln f_2)\,g(U, V)$. Combining this with
 \eqref{4.13} gives $\sin^2\theta\,X(\ln f_2)\,g(U, V)=0$.
Because $g$ is Riemannian metric,  this implies that $\sin^2\theta\, X(\ln f_2)=0.$ Since $N_\theta$ is proper pointwise slant, we must have
$X(\ln f_2)=0$. Thus $f_2$ is constant. 

Conversely, if $f_2$ is constant, then  \eqref{4.2} gives
 $g(h(\mathfrak{D}^{T}, \mathfrak{D}^\theta),F\mathfrak{D}^\theta)=0$.
Further, it follows from the second equality of Lemma \ref{L:4.3} that
$g(h(\mathfrak{D}^{T}, \mathfrak{D}^\theta), \varphi \mathfrak{D}^\perp)=0.$
Moreover, we find from the hypothesis that
$g(h(\mathfrak{D}^{T}, \mathfrak{D}^\theta), \mu)=0$. Consequently, $M^n$ is ${\mathfrak{D}^{T}}\oplus\mathfrak{D}^\theta$-mixed totally geodesic.
\end{proof}


\section{A general inequality for bi-warped products}\label{S5}

Let $M^n=N_T\times_{f_1}N^{n_{1}}_\perp\times_{f_2}N^{n_{2}}_\theta$ be a $CRS$ bi-warped product in a cosymplectic manifold $\widetilde M^{2m+1}$ with $\dim N_{T}=2p+1$. 
We choose orthonormal frames of ${\mathfrak{D}^{T}},\, {\mathfrak{D}}^\perp$ and ${\mathfrak{D}}^\theta$ of $TM^n$ given respectively by
\begin{align*}
&{\mathfrak{D}^{T}}=\rm{Span}\{e_1,\,\cdots,e_{p},\,e_{p+1}=\varphi e_{1},\,\cdots,e_{2p}=\varphi e_{p}, e_{2p+1}={\mathfrak{\xi}}\},\;\\&{\mathfrak{D}}^\perp=\rm{Span}\{e_{2p+2}=\bar e_1,\,\cdots,\,e_{2p+1+n_1}=\bar e_{n_1}\},\\
&{\mathfrak{D}}^\theta=\rm{Span}\{e_{2p+2+n_1}=e_1^{*},\,\cdots,\,e_{2p+1+n_1+s}=e_s^{*},\,\\&\hskip.5in e_{2p+2+n_1+s}=e_{s+1}^{*}=\sec\theta\, Te_1^{*},\cdots, e_{n}=e_{2s}^{*}=\sec\theta\, Te_{s}^{*}\}.
\end{align*}
Then we have the following orthonormal frames for $\varphi{\mathfrak{D}}^\perp,\,F\mathfrak{D}^\theta$ and $\mu$ as follows:   
\begin{align*}
&\varphi{\mathfrak{D}}^\perp=\rm{Span}\{e_{n+1}= \varphi\bar e_1,\cdots,\, e_{n+n_1}=\varphi\bar e_{n_1}\},\\
&F\mathfrak{D}^\theta=\rm{Span}\{e_{n+n_1+1}=\hat{e}_{1}=\csc\theta\, Fe_{1}^{*},\cdots,\,e_{n+n_1+s}=\hat{e}_{s}=\csc\theta\, Fe_{s}^{*},\\
&\hskip.2in e_{n+n_1+s+1}=\hat{e}_{s+1}=\csc\theta\sec\theta\, FTe_{1}^{*},\cdots,\, e_{n+n_1+n_2}=\hat{e}_{n_2}=\csc\theta\sec\theta\, FTe_{s}^{*}\},\\
&\mu=\rm{Span}\{e_{n+n_1+n_2+1}=\widetilde e_{1},\cdots, e_{2m+1}=\widetilde e_{2m+1-n-n_1-n_2}\}.
\end{align*}
Clearly,  ${\rm rank}\;{\varphi\mathfrak{D}}^\perp=n_1,\,\,{\rm rank}\; F{\mathfrak{D}}^\theta=2s=n_2$ and ${\rm rank}\;\mu=2m+1-n-n_1-n_2$.

The following is the main result of this section.

\begin{theorem}\label{T:5.1} Let $M=N_T\times_{f_1}N^{n_{1}}_\perp\times_{f_2}N^{n_{2}}_\theta$ be a $CRS$ bi-warped product in a $(2m+1)$-dimensional cosymplectic manifold $\widetilde M^{2m+1}$. Then
\begin{enumerate}
\item[(i)] The squared norm of the second fundamental from satisfies 
\begin{equation}\label{5.1}
\|h\|^2\geq 2n_1\|\nabla(\ln f_1)\|^2+2n_2(1+2\cot^2\theta)\,\|\nabla(\ln f_2)\|^2, 
\end{equation}
where $\nabla(\ln f_1)$ and $\nabla(\ln f_2)$ are the gradients of $\ln f_{1}$ and $\ln f_{2}$, respectively.
\item[(ii)] If the equality sign holds in \eqref{5.1}, then $N_T$ is totally geodesic in $\widetilde M^{2m+1}$ and $N^{n_{1}}_\perp$ and $N^{n_{2}}_\theta$ are totally umbilical submanifolds of $\widetilde M^{2m+1}$. Moreover, $M$ is a ${\mathfrak{D}}^\perp\oplus{\mathfrak{D}}^\theta$-mixed totally geodesic submanifold but neither ${\mathfrak{D}^{T}}\oplus {\mathfrak{D}}^\perp$-mixed totally geodesic nor ${\mathfrak{D}^{T}}\oplus {\mathfrak{D}}^\theta$-mixed totally geodesic.
\end{enumerate}
\end{theorem}
\begin{proof} From \eqref{funda}, we have
\begin{align*}
\|h\|^2=\sum_{i,j=1}^ng(h(e_i, e_j), h(e_i, e_j)) =\sum_{r=n+1}^{2m+1}\sum_{i,j=1}^ng(h(e_i, e_j), e_r)^2.\end{align*}
By using \eqref{tan-nor}, the above equation takes the form
\begin{align*}
\|h\|^2&=\sum_{r=n+1}^{2m+1}\sum_{i,j=1}^{2p+1}g(h(e_i, e_j), e_r)^2+\sum_{r=n+1}^{2m+1}\sum_{i,j=1}^{n_1}g(h(\bar e_i, \bar e_j), e_r)^2\\
&+\sum_{r=n+1}^{2m+1}\sum_{i,j=1}^{n_2}g(h(e_i^{*}, e_j^{*}), e_r)^2+2\sum_{r=n+1}^{2m+1}\sum_{i=1}^{2p+1}\sum_{j=1}^{n_1}g(h(e_i, \bar e_j), e_r)^2\\
&+2\sum_{r=n+1}^{2m+1}\sum_{i=1}^{2p+1}\sum_{j=1}^{n_2}g(h(e_i, e_j^{*}), e_r)^2+2\sum_{r=n+1}^{2m+1}\sum_{i=1}^{n_1}\sum_{j=1}^{n_2}g(h(\bar e_i,  e_j^{*}), e_r)^2.
\end{align*}
Again, we derive from \eqref{tan-nor} that
\begin{align}
\|h\|^2&=\sum_{r=1}^{n_1}\sum_{i,j=1}^{2p+1}g(h(e_i, e_j), \varphi \bar e_r)^2+\sum_{r=1}^{n_2}\sum_{i,j=1}^{2p+1}g(h(e_i, e_j), \hat{e}_r)^2\notag\\
&+\sum_{r=1}^{2m+1-n-n_1-n_2}\sum_{i,j=1}^{2p+1}g(h(e_i, e_j), \widetilde e_r)^2+\sum_{r=1}^{n_1}\sum_{i,j=1}^{n_1}g(h(\bar e_i, \bar e_j), \varphi \bar e_r)^2
\notag\\
&+\sum_{r=1}^{n_2}\sum_{i,j=1}^{n_1}g(h(\bar e_i, \bar e_j), \hat{e}_r)^2+\sum_{r=1}^{2m+1-n-n_1-n_2}\sum_{i,j=1}^{n_1}g(h(\bar e_i, \bar e_j), \widetilde {e}_r)^2\notag\\
&+\sum_{r=1}^{n_1}\sum_{i,j=1}^{n_2}g(h(e_i^{*}, e_j^{*}), \varphi \bar e_r)^2+\sum_{r=1}^{n_2}\sum_{i,j=1}^{n_2}g(h(e_i^{*}, e_j^{*}), \hat{e}_r)^2
\notag\\ \label{5.2}
&+2\sum_{r=1}^{n_1}\sum_{i=1}^{2p+1}\sum_{j=1}^{n_1}g(h(e_i, \bar e_j), \varphi \bar e_r)^2+2\sum_{r=1}^{n_2}\sum_{i=1}^{2p+1}\sum_{j=1}^{n_1}g(h(e_i, \bar e_j), \hat e_r)^2 
\\
&+2\sum_{r=1}^{2m+1-n-n_1-n_2}\sum_{i=1}^{2p+1}\sum_{j=1}^{n_1}g(h(e_i, \bar e_j), \widetilde e_r)^2+2\sum_{r=1}^{n_2}\sum_{i=1}^{2p+1}\sum_{j=1}^{n_2}g(h(e_i, e_j^{*}), \hat e_r)^2\notag
\\& \notag
+2\sum_{r=1}^{n_1}\sum_{i=1}^{2p+1}\sum_{j=1}^{n_2}g(h(e_i, e_j^{*}), \varphi\bar e_r)^2+2\sum_{r=1}^{2m+1-n-n_1-n_2}\sum_{i=1}^{2p+1}\sum_{j=1}^{n_2}g(h(e_i, e_j^{*}), \widetilde e_r)^2\notag\\
&+2\sum_{r=1}^{n_1}\sum_{i=1}^{n_1}\sum_{j=1}^{n_2}g(h(\bar e_i,  e_j^{*}), \varphi \bar e_r)^2+2\sum_{r=1}^{n_2}\sum_{i=1}^{n_1}\sum_{j=1}^{n_2}g(h(\bar e_i,  e_j^{*}), \hat e_r)^2
\notag\\
&+2\sum_{r=1}^{2m+1-n-n_1-n_2}\sum_{i=1}^{n_1}\sum_{j=1}^{n_2}g(h(\bar e_i,  e_j^{*}), \widetilde e_r)^2.\notag
\end{align}
Let us omit the  positive 3rd, 6th, 9th, 12th, 15th and 18th terms of $\mu$-components in the right hand side of \eqref{5.2}. Since we couldn't find any relations for the 4th, 5th, 7th, 8th, 16th and 17th terms, we also omit these positive terms and using Lemma \ref{L:4.2}(i), the 1st and 2nd terms vanish identically. Similarly, the 11th and 13th terms vanish according to Lemma \ref{L:4.3}. Hence, with the remaining 10th and 14th terms, the above expression becomes
\begin{equation}\begin{aligned}\label{5.3}
\|h\|^2\geq\; & 2\sum_{r=1}^{n_1}\sum_{i=1}^{p}\sum_{j=1}^{n_1}g(h(e_i, \bar e_j), \varphi \bar e_r)^2+2\sum_{r=1}^{n_1}\sum_{i=1}^{p}\sum_{j=1}^{n_1}g(h(\varphi e_i, \bar e_j), \varphi \bar e_r)^2\\ 
&+2\sum_{r=1}^{n_1}\sum_{j=1}^{n_1}g(h(\mathfrak{\xi}, \bar e_j), \varphi \bar e_r)^2+2\sum_{r=1}^{n_2}\sum_{i=1}^{p}\sum_{j=1}^{n_2}g(h(e_i, e_j^{*}), \hat e_r)^2
\\& +2\sum_{r=1}^{n_2}\sum_{i=1}^{p}\sum_{j=1}^{n_2}g(h(\varphi e_i, e_j^{*}), \hat e_r)^2+2\sum_{r=1}^{n_2}\sum_{j=1}^{n_2}g(h((\mathfrak{\xi}, e_j^{*}), \hat e_r)^2.
\end{aligned}\end{equation}

From Lemma \ref{L:4.1}, the 3rd and 6th terms of \e{5.3} are identical zero. Then, by using (ii)-(iii) of Lemma \ref{L:4.2} and  \eqref{4.2}-\eqref{4.8} with orthonormality of vector fields, we find
\begin{align}\label{5.4}
\|h\|^2\; & \geq2n_1\sum_{i=1}^{p}(e_i\ln f_1)^2+2n_1\sum_{i=1}^{p}(\varphi e_i\ln f_1)^2+2n_2(\csc^2\theta+\cot^2\theta)\sum_{i=1}^{p}(\varphi e_i\ln f_2)^2\notag\\
&+2n_2(\csc^2\theta+\cot^2\theta)\sum_{i=1}^{p}(e_i\ln f_2)^2.
\end{align}
Since  $\mathfrak{\xi}(\ln f_1)=\mathfrak{\xi}(\ln f_2)=0$ from Lemma \ref{L:4.1}, the above relation becomes
\begin{equation}\begin{aligned}\label{5.5}
\|h\|^2&\geq2n_1\sum_{i=1}^{2p+1}(e_i\ln f_1)^2+2n_2(\csc^2\theta+\cot^2\theta)\sum_{i=1}^{2p+1}(e_i\ln f_2).
\end{aligned}\end{equation}
So,  we find the inequality given in statement (i) via the definition of gradient. 

For the equality, it follows from the omitted $\mu$-components terms in \eqref{5.2} that\begin{align}\label{eq1}
h(TM, TM)\perp\mu
\end{align}
Also, from the omitted 4th and 5th terms in \eqref{5.2}, we have
\begin{align}\label{eq2}
h({\mathfrak{D}}^\perp, {\mathfrak{D}}^\perp)\perp \varphi{\mathfrak{D}}^\perp,\,\,\,\,\,h({\mathfrak{D}}^\perp, {\mathfrak{D}}^\perp)\perp F{\mathfrak{D}}^\theta.
\end{align}
Then  we obtain from \eqref{eq1} and \eqref{eq2} that
\begin{align}\label{eq3}
h({\mathfrak{D}}^\perp, {\mathfrak{D}}^\perp)=\{0\}.
\end{align}
From the omitted 7th and 8th terms in the right hand side of \eqref{5.2}, we get
\begin{align}\label{eq5}
h({\mathfrak{D}}^\theta, {\mathfrak{D}}^\theta)\perp \varphi{\mathfrak{D}}^\perp,\,\,\,\,\,h({\mathfrak{D}}^\theta, {\mathfrak{D}}^\theta)\perp F{\mathfrak{D}}^\theta.
\end{align}
Thus  we conclude from \eqref{eq1} and \eqref{eq5} that
\begin{align}\label{eq6}
h({\mathfrak{D}}^\theta, {\mathfrak{D}}^\theta)=\{0\}.
\end{align}

Similarly, from the omitted 16th and 17th terms in the right hand side of  \eqref{5.2}, we arrive at
\begin{align}\label{eq7}
h({\mathfrak{D}}^\perp, {\mathfrak{D}}^\theta)\perp\varphi{\mathfrak{D}}^\perp,\,\,\,\,\,h({\mathfrak{D}}^\perp, {\mathfrak{D}}^\theta)\perp F{\mathfrak{D}}^\theta.
\end{align}
From \eqref{eq1} and \eqref{eq7}, we find
\begin{align}\label{eq8}
h({\mathfrak{D}}^\perp, {\mathfrak{D}}^\theta)=\{0\}.
\end{align}

On the other hand, we find from the vanishing 1st and 2nd terms of \e{5.2} that
\begin{align}\label{eq9}
h({\mathfrak{D}^{T}}, {\mathfrak{D}^{T}})\perp\varphi{\mathfrak{D}}^\perp,\,\,\,\,\,h({\mathfrak{D}^{T}},\, {\mathfrak{D}^{T}})\perp F{\mathfrak{D}}^\theta.
\end{align}
Thus,  we derive from \eqref{eq1} and \eqref{eq9}, that
\begin{align}\label{eq10}
h({\mathfrak{D}^{T}}, {\mathfrak{D}^{T}})=\{0\}.
\end{align}
And, from the vanishing 11th term of \eqref{5.2} and \eqref{eq1}, we get
\begin{align}\label{eq11}
h({\mathfrak{D}^{T}}, {\mathfrak{D}}^\perp)\subset \varphi{\mathfrak{D}}^\perp.
\end{align}
Similarly, from the vanishing 13th term in \eqref{5.2} with \eqref{eq1}, we find
\begin{align}\label{eq12}
h({\mathfrak{D}^{T}}, {\mathfrak{D}}^\theta)\subset F{\mathfrak{D}}^\theta.
\end{align}
Since $N_T$ is totally geodesic in $M$ (see Remark \ref{R:3.1}), using this fact with \eqref{eq3}, \eqref{eq6} and \eqref{eq10}. Again, from \eqref{eq8}, \eqref{eq11} and \eqref{eq12} with Remark \ref{R:3.1}, we conclude that $N_\perp$ and $N_\theta$ are totally umbilical in $\widetilde M$, while; using all conditions with  \eqref{eq8} $M$ is a ${\mathfrak{D}}^\perp\oplus {\mathfrak{D}}^\theta$-mixed totally geodesic submanifold of $\widetilde M$ but not ${\mathfrak{D}^{T}}\oplus {\mathfrak{D}}^\perp$ and ${\mathfrak{D}^{T}}\oplus {\mathfrak{D}}^\theta$-mixed totally geodesic. 
\end{proof}

If $n_{2}=0$, then  $N_T\times_{f_1}N^{n_{1}}_\perp\times_{f_2}N^{n_{2}}_\theta$ reduces to a contact CR-warped product $N_T\times_{f_1}N^{n_{1}}_\perp$. Thus,
Theorem \ref{T:5.1} implies to the following main result of \cite{UA}.

\begin{corollary}\label{C:5.1}  Let $M=N_T\times_{f}N^{n_{1}}_\perp$ be a contact CR-warped product of a cosymplectic manifold. Then 
$\|h\|^2\geq2n_1\! \|\nabla(\ln f)\|^2$.
\end{corollary}

If $n_{1}=0$, then $N_T\times_{f_1}N^{n_{1}}_\perp\times_{f_2}N^{n_{2}}_\theta$ reduces to a pointwise semi-slant warped product $N_T\times_{f_{2}}N_\theta$. In this case,
Theorem \ref{T:5.1} implies to the main result of \cite{Park}.

\begin{corollary}\label{C:5.2}  Let $M=N_T\times_{f}N^{n_{2}}_\theta$ be a  pointwise semi-slant warped product in a cosymplectic manifold. Then 
$\|h\|^2\geq2n_2(\csc^2\theta+\cot^2\theta)\|\nabla(\ln f)\|^2.$
\end{corollary}

\begin{remark} By applying the same methods given in \cite{CD} and in \cite{HM}, we have the following extension of Corollary \ref{C:5.1} to  contact multiply CR-warped product  $N_T\times_{f_1}N_1^{n_{1}}\times\cdots\times_{f_k}N_k^{n_{k}}$ of cosymplectic manifolds $\widetilde M^{2m+1}$, where $N_T$ is invariant and $N_{1}^{1},\ldots,N_{k}^{n_{k}}$ are anti-invariant submanifolds of $\widetilde M^{2m+1}$.

\begin{theorem}\label{T:5.2} Let $M=N_T\times_{f_1}N_1^{n_{1}}\times\cdots\times_{f_k}N_k^{n_{k}}$ be a contact multiply CR-warped product in a cosymplectic manifold $\widetilde M^{2m+1}$. Then
\begin{align}\label{5.17}
\|h\|^2\geq2\sum_{i=1}^{k}n_i \|\nabla(\ln f_i)\|^2.
\end{align}

The equality sign in \eqref{5.17} holds identically if and only if the following three conditions hold:
\begin{enumerate}
\item[(a)] $N_T$ is a totally geodesic submanifold of $\widetilde M^{2m+1}$.
\item[(b)] $N_i\; (i=1,\cdots, k)$ are totally umbilical submanifolds of $\widetilde M^{2m+1}$.
\item[(c)] $_{f_1}N_1\times\cdots\times_{f_k}N_k$ is mixed totally geodesic in $\widetilde M^{2m+1}$
\end{enumerate}
\end{theorem}
\end{remark}

\section{Further applications of Theorem \ref{T:5.1}}\label{S6}

The Dirichlet energy of a function $\psi$ on a compact manifold $M$ is defined by
\begin{align}\label{eng}E(\psi)=\frac{1}{2}\int_{M}\|\nabla\psi\|^2\,dV, \end{align}
where $\nabla\psi$ is the gradient of $\psi$ and $dV$ is the volume element.

Further applications of Theorem \ref{T:5.1} is to express the Dirichlet energy of the warping functions $f_1$ and $f_2$. In order to do so, we assume that the invariant factor $N_{T}$ of the $CRS$ bi-warped product is compact without boundary.

\begin{theorem}\label{T:6.1} Let $M^n=N_T\times_{f_1}N^{n_1}_\perp\times_{f_2}N^{n_2}_\theta$ be a $CRS$ bi-warped product in a cosymplectic manifold $\widetilde M^{2m+1}$ with compact $N_T$.  Then,  for each $p\in N^{n_{1}}_\perp$ and $q\in N^{n_{2}}_\theta$, we have the following inequality for the Dirichlet energy of $f_1$ and $f_2$:
\begin{equation}\begin{aligned}\notag 
n_1E(\ln f_1)+n_2(1+2\cot^{2}\theta)E(\ln f_2) \leq\frac{1}{4}\int_{N_T\times\{p\}\times\{q\}}\|h\|^2\,dV_T.
\end{aligned}\end{equation}
 \end{theorem}
\begin{proof} By integrating \eqref{5.1} over $M$, we derive the required inequality. 
\end{proof}

The following two corollaries are immediate consequences of Theorem \ref{T:6.1}. 

\begin{corollary}\label{C:6.1} Let $M=N_T\times_fN^{n_{1}}_\perp$ be a contact CR-warped product of a cosymplectic manifold $\widetilde M^{2m+1}$. If $N_T$ is compact and  $q\in N^{n_{1}}_\perp$, then
\begin{align}\notag
E(\ln f_1)\le\frac{1}{4n_1}\int_{N_T\times\{q\}}\|h\|^2\,dV_T.
\end{align}
\end{corollary}

Similarly, if $\dim N^{n_{1}}_\perp=0$, then Theorem \ref{T:6.1} implies the following.

\begin{corollary}\label{C:6.2} Let $M=N_T\times_fN^{n_{2}}_\theta$ be a  pointwise semi-slant warped product in a cosymplectic manifold $\widetilde M^{2m+1}$. If $N_T$ is compact and $s\in N^{n_{2}}_\theta$, then
\begin{align}\notag
E(\ln f_2)\le\frac{1}{4n_2(1+2\cot^2\theta)}\int_{N_T\times\{s\}}\|h\|^2\,dV_T.
\end{align}
\end{corollary}

\section{A non-trivial example of $CRS$-warped product submanifolds}\label{S7}

In this section, we provide a non-trivial example of $CRS$ bi-warped products in a cosymplectic manifold satisfying the equality case of inequality \eqref{5.1} identically.
\begin{example}
\rm{Let ${\mathbb{R}}^{19}$ denotes the Cartesian space with the Cartesian coordinates $(x_1,\,y_1,\,\cdots,\,x_9,\, y_9,\, z)$. Consider the contact metric structure $(\varphi,\xi,\eta,g)$  on ${\mathbb{R}}^{19}$, where $g$ is the standard Euclidean metric and $\varphi,\xi,\eta$ are defined by
\begin{equation}\label{7.1}
\xi=\frac{\partial}{\partial z},\; \eta=dz,\;  \varphi\left(\frac{\partial}{\partial x_i}\right)=-\frac{\partial}{\partial y_i},\;  \varphi\left(\frac{\partial}{\partial y_i}\right)=\frac{\partial}{\partial x_i},\; \varphi \left(\frac{\partial}{\partial z}\right)=0, 
\end{equation} 
for $i=1,\ldots,9$. Also,  we have $d\eta=d\varphi =0$ and $\widetilde \nabla \varphi=\widetilde \nabla \xi=0$. Hence, ${\mathbb{R}}^{19}$ is a cosymplectic manifold.

Let $M^{7}$ be the submanifold of ${\mathbb{R}}^{19}$ given by 
\begin{equation}\begin{aligned}\label{7.2}
\psi(u,\, v,\, w,\, s,\, t,\, r,\, z)=\; &\Big(u\cos w,\, v\cos w,\, u\cos s,\, v\cos s,\,u\sin w,\, v\sin w,\, \\
&\hskip-1.1in  u\sin s,\, v\sin s, u\cos t,\, v\cos t,\, u\cos r,\,v\cos r,\, u\sin t,\, v\sin t,\, u\sin r,\,\\
&\hskip-.7in   v\sin r,\,k(r-t),\, -k(r+t), z\Big), \;\; u\neq 0,\; v\neq 0,\; k\in\mathbb{R}-\left\lbrace 0 \right\rbrace.
\end{aligned}\end{equation}
Then the tangent bundle $TM^{7}$ is spanned by the following coordinate vector fields with respect to $u,\, v,\, w,\, s,\, t,\, r,\, z$, respectively.
\begin{align}\notag
&W_1=\cos w\,\frac{\partial}{\partial x_1}+\cos s\,\frac{\partial}{\partial x_2}+\sin w\,\frac{\partial}{\partial x_3}+\sin s\,\frac{\partial}{\partial x_4}+\cos t\,\frac{\partial}{\partial x_5}+\cos r\,\frac{\partial}{\partial x_6}
\\&\notag\hskip.3in +\sin t\,\frac{\partial}{\partial x_7}+\sin r\,\frac{\partial}{\partial x_8},\\
&\notag W_2=\cos w\,\frac{\partial}{\partial y_1}+\cos s\,\frac{\partial}{\partial y_2}+\sin w\,\frac{\partial}{\partial y_3}+\sin s\,\frac{\partial}{\partial y_4}+\cos t\,\frac{\partial}{\partial y_5} +\cos r\,\frac{\partial}{\partial y_6}\\&\notag\hskip.3in+\sin t\,\frac{\partial}{\partial y_7}+\sin r\,\frac{\partial}{\partial y_8},\\
&\label{7.3} W_3=-u\sin w\,\frac{\partial}{\partial x_1}-v\sin w\,\frac{\partial}{\partial y_1}+u\cos w\,\frac{\partial}{\partial x_3}+v\cos w\,\frac{\partial}{\partial y_3},\\
 &\notag W_4=-u\sin s\,\frac{\partial}{\partial x_2}-v\sin s\,\frac{\partial}{\partial y_2}+u\cos s\,\frac{\partial}{\partial x_4}+v\cos s\,\frac{\partial}{\partial y_4},\\
&\notag W_5=-u\sin t\,\frac{\partial}{\partial x_5}-v\sin t\,\frac{\partial}{\partial y_5}+u\cos t\,\frac{\partial}{\partial x_7}+v\cos t\,\frac{\partial}{\partial y_7}
 -k\frac{\partial}{\partial x_9}-k\frac{\partial}{\partial y_9},\\
&\notag W_6=-u\sin r\,\frac{\partial}{\partial x_6}-v\sin r\,\frac{\partial}{\partial y_6}+u\cos r\,\frac{\partial}{\partial x_8}+v\cos r\,\frac{\partial}{\partial y_8}
+k\frac{\partial}{\partial x_9}-k\frac{\partial}{\partial y_9},
\\&\notag W_7=\frac{\partial}{\partial z}.
\end{align}
Then we find $\varphi W_7=0$ and
\begin{align}\notag
&\varphi W_1=-\cos w\,\frac{\partial}{\partial y_1}-\cos s\,\frac{\partial}{\partial y_2}-\sin w\,\frac{\partial}{\partial y_3}-\sin s\,\frac{\partial}{\partial y_4}-\cos t\,\frac{\partial}{\partial y_5}
\\&\notag \hskip.45in -\cos r\,\frac{\partial}{\partial y_6}-\sin t\,\frac{\partial}{\partial y_7}-\sin r\,\frac{\partial}{\partial y_8},
\\&\notag  \varphi W_2=\cos w\,\frac{\partial}{\partial x_1}+\cos s\,\frac{\partial}{\partial x_2}+\sin w\,\frac{\partial}{\partial x_3}+\sin s\,\frac{\partial}{\partial x_4}+\cos t\,\frac{\partial}{\partial x_5}
\\&\notag \hskip.45in +\cos r\,\frac{\partial}{\partial x_6}+\sin t\,\frac{\partial}{\partial x_7}+\sin r\,\frac{\partial}{\partial x_8},
 \\&\label{7.4} \varphi W_3=-v\sin w\,\frac{\partial}{\partial x_1}+u\sin w\,\frac{\partial}{\partial y_1}+v\cos w\,\frac{\partial}{\partial x_3}-u\cos w\,\frac{\partial}{\partial y_3},
\\ &\notag  \varphi W_4=-v\sin s\,\frac{\partial}{\partial x_2}+ u\sin s\,\frac{\partial}{\partial y_2}+v\cos s\,\frac{\partial}{\partial x_4}-u\cos s\,\frac{\partial}{\partial y_4},
\\&\notag \varphi W_5=-v\sin t\frac{\partial}{\partial x_5}+u\sin t\,\frac{\partial}{\partial y_5}+v\cos t\frac{\partial}{\partial x_7}-u\cos t\frac{\partial}{\partial y_7}-k\frac{\partial}{\partial x_9}
 +k\frac{\partial}{\partial y_9}, \\
&\notag \varphi W_6=-v\sin r\,\frac{\partial}{\partial x_6}+u\sin r\,\frac{\partial}{\partial y_6}+v\cos r\,\frac{\partial}{\partial x_8}-u\cos r\,\frac{\partial}{\partial y_8}-k\frac{\partial}{\partial x_9}  -k\frac{\partial}{\partial y_9}.
\end{align}
Clearly, ${\mathfrak{D}^{T}}=\rm{Span}\{W_1,W_2,W_{7}\}$ is an invariant distribution. Since both $\varphi W_3$ and $\varphi W_4$ are orthogonal to $TM^{7}$, ${\mathfrak{D}}^{\perp}=\rm{Span}\{W_3, W_4\}$ is an anti-invariant distribution. Further, it is direct to show that ${\mathfrak{D}^{\theta}}=\rm{Span}\{W_5,W_6\}$ is a pointwise slant distribution with slant function $\theta=\arccos (\frac{2k^2}{u^2+v^2+2k^2})$. 

Obviously, ${\mathfrak{D}^{T}},\, {\mathfrak{D}}^{\perp}$ and ${\mathfrak{D}}^{\theta}$ are integrable distributions on $M$.  If we denote the integral manifolds of ${\mathfrak{D}^{T}},\;{\mathfrak{D}}^\perp$ and ${\mathfrak{D}^{T}}^{\theta}$ by $N_T,\, N_\perp$ and $N_\theta$, respectively, then the induced metric tensor $g$ of $M= N_T\times_{f_1}N_\perp\times_{f_2} N_\theta$ is given by
\begin{equation}\begin{aligned}\label{7.5}
&g=g_{N_T}+f^2_1g_{N_\perp}+f^2_2g_{N_\theta},\\
&g_{N_T}=4(du^2+dv^2)+dz^2,\;\; g_{N_{\perp}}=dw^2+ds^2,\;\; g_{N_{\theta}}=dt^2+dr^2,
\\& f_{1}=\sqrt{u^{2}+ v^{2}}, \;\; f_2=\sqrt{u^2+v^2+2k^2}.
\end{aligned}\end{equation}
Hence, $M^{7}$ is a $CRS$ bi-warped product in ${\mathbb{R}}^{19}$ which satisfy $\xi f_{1}=\xi f_{2}=0$. 

Let $\widetilde \nabla$ denote the Riemannian connection on ${\mathbb{R}}^{19}$.
Since $h(W_{i},W_{j})$ is the normal component of $\widetilde \nabla_{W_{i}}W_{j}$, it follows easily from \e{7.3} that
\begin{align}\label{7.6} &h(\mathfrak{D}^{T}, \mathfrak{D}^{T})=h(\mathfrak{D}^{\perp}, \mathfrak{D}^{\perp})=h(\mathfrak{D}^{\perp}, \mathfrak{D}^{\theta})=h(\mathfrak{D}^{\theta}, \mathfrak{D}^{\theta})=\{0\},
\\&\label{7.6}h(\xi, \mathfrak{D}^{\perp})=h(\xi, \mathfrak{D}^{\theta})=\{0\},\quad \xi=\frac{\partial}{\partial z}.
\end{align}

  Furthermore, we find from \e{7.3} and \e{7.4} that 
\begin{align*}
&\widetilde \nabla_{W_{1}}W_{3}=\frac{u}{u^{2}+ v^{2}}W_{3}+\frac{uv}{(u^{2}+ v^{2})^{2}}\varphi W_{3},\end{align*}
which implies $0\ne h(W_{1},W_{3})\in \varphi(\mathfrak{D}^{\perp})$. In the same way, we  find
\begin{align}\label{7.8}
h({\mathfrak D}^{T}, \mathfrak{D}^\perp)\neq\{0\}\;\; {\rm and}\;\; h(\mathfrak{D}^{T}, \mathfrak{D}^{\perp})\in \varphi(\mathfrak{D}^{\perp}).
\end{align}
By computing  $\widetilde \nabla_{W_{i}}W_{a}$ for $i=1,2$ and $a=5,6$, and by applying \e{7.6}, we also have
\begin{align}\label{7.9}
h({\mathfrak D}^{T}, \mathfrak{D}^\theta)\neq\{0\}\;\; {\rm and}\;\; h(\mathfrak{D}^{T}, \mathfrak{D}^{\theta})\in F(\mathfrak{D}^{\theta}).
\end{align}
Consequently, Remark \ref{R:3.1} and \eqref{7.6}--\eqref{7.9} imply that $N_T$ is totally geodesic in $\mathbb{R}^{19}$ and $M$ is $\mathfrak{D}^\perp\oplus\mathfrak{D}^\theta$-mixed totally geodesic. Furthermore, it follows from \eqref{7.8} and \eqref{7.9} that $M^{7}$ is neither $\mathfrak{D}^{T}\oplus\mathfrak{D}^\perp$-mixed totally geodesic nor $\mathfrak{D}^{T}\oplus\mathfrak{D}^\theta$-mixed totally geodesic in ${\mathbb{R}}^{19}$.}
\end{example}

\vskip.1in

\noindent {\bf Acknowledgement.} The authors thank Professor Kwang Soon Park for pointing out an error in an earlier version of this article.

\end{document}